# Krasovskiĭ Stability Theorem for FDEs in the Extended Sense

Qian Feng[1] and Wilfrid Perruquetti[2]

***Abstract*— The analysis of the stability of systems' equilibria plays a central role in the study of dynamical systems and control theory. This note establishes an extension of the celebrated Krasovskiĭ stability theorem for functional differential equations (FDEs) in the extended sense. Namely, the FDEs hold for $t \geq t_0$ almost everywhere with respect to the Lebesgue measure. The existence and uniqueness of such FDEs were briefly discussed in J.K Hale's classical treatise on FDEs, yet a corresponding stability theorem was not provided. A key step in proving the proposed stability theorem was to utilize an alternative strategy instead of relying on the mean value theorem of differentiable functions. The proposed theorem can be useful in the stability analysis of cybernetic systems, which are often subject to noise and glitches that have a countably infinite number of jumps. To demonstrate the usefulness of the proposed theorem, we provide examples of linear systems with time-varying delays in which the FDEs cannot be defined in the conventional sense.**

## I. INTRODUCTION

Stability is crucial for qualitatively describing the asymptotic behavior of system dynamics, as explicit solutions to the equations of nonlinear systems cannot usually be derived. As pioneered in the seminal works of Aleksandr Lyapunov [1]–[3], the stability of the equilibria of ODEs can be determined by the existence of an energy-like positive definite function, whose time derivative along the system's trajectory is negative definite. However, finding such Lyapunov functions is often challenging, and only sufficient conditions can be formulated for their existence. Nevertheless, the direct method of Lyapunov has become an indispensable tool [4], [5] for the analysis and control of nonlinear systems.

The original idea of the Lyapunov direct method was later extended by Nikolaĭ Nikolayevich Krasovskiĭ [6], [7] to address the stability of systems with delayed arguments, and later to functional differential equations [8]–[11]. As the state space of functional differential equations is infinite-dimensional, the Lyapunov functions in this case are replaced by the Krasovskiĭ functionals, whose construction inevitably involves the use of infinite-dimensional analysis. This explains why, even for linear time-delay systems [12], the construction of Krasovskiĭ functionals remains challenging, in contrast to linear ODEs with finite dimensions whose stability can be addressed by constructing a quadratic Lyapunov function.

Similar to Lyapunov's achievements, the construction of Krasovskiĭ functionals has become a standard approach for the analysis of time-delay systems in the control community. The Krasovskiĭ stability theorem [13, Section 5.2] can be applied to FDEs satisfying a list of prerequisites, including the continuity (piecewise) with respect to both the time and state arguments of the right-hand side of the FDEs to ensure that the resulting solution $\boldsymbol{x}(t)$ is differentiable. Some of the aforementioned prerequisites, however, can be conservative in applied settings, as many engineering systems are often subject to noise and glitches which can lead to a continuous yet not strictly differentiable solution $\boldsymbol{x}(t)$. An illustrative example is when systems have a time-varying delay $r(t)$ that is integrable yet **not piecewise continuous**, which cannot be addressed by the Krasovskii stability theorems formulated for time-varying delay systems in [14], [15].

The above problem can be addressed using the Carathéodory formulation of FDEs [13, Theorem 5.3, Chapter 1], similar to the case of ODEs [16]. Namely, the derivative in the FDEs can be interpreted as weak derivative, which holds for $t \geq t_0$ almost everywhere with respect to the Lebesgue measure. This means that a solution to the FDEs only needs to be locally absolutely continuous, which is differentiable almost everywhere. To ensure the existence and uniqueness of the initial value problem for the FDEs, the right-hand side of the FDE must satisfy the Carathéodory conditions [13, section 2.6] and is locally Lipschitz [17] in the state argument. As a result, continuity in the time argument is no longer required, but it needs to be measurable instead. The above characterization of FDEs can be interpreted as an extended version of the conventional framework using ordinary derivatives. A countable number of discontinuities in the right-hand side of FDEs can be handled by the weak derivative framework, which is a notable advantage in the modeling of engineering and cybernetics systems.

The Krasovskiĭ stability theorem in most literature [13, Section 5.2], [18, Theorem 1.3], [19, Theorem 1.1], [11, Theorem 3.9], [20, Theorem 7.2.2.], [Section 4.2] [21], [22, Theorem 5.1], however, is formulated using the conventional derivatives, which means that solutions must be differentiable for all $t \geq t_0$ (at least right-hand differentiable) and the right-hand side of the FDEs must be piecewise continuous with respect to time. Given the usefulness of the Carathéodory formulation of FDEs, it would be desirable to have a corresponding Krasovskiĭ stability theorem, allowing the stability of the extended version of FDEs to be analyzed in a similar manner.

This work was partially supported by the National Natural Science Foundation of China under Grant Nos. 62303180 and 62273145, and Fundamental Research Funds for Central Universities under Grant 2023MS032, China, and ANR (France) Project ANR-15-CE23-0007.

1. School of Control and Computer Engineering, North China Electric Power University Beijing, China. Email: `qianfeng@ncepu.edu.cn`, `qfen204@aucklanduni.ac.nz`

2. École Centrale de Lille, Cité Scientifique, 59650 Villeneuve-d'Ascq, Lille, France Emails: `wilfrid.perruquetti@centralelille.fr`

In this short note, we formally establish the Krasovskiĭ stability theorem for FDEs in the extended sense, where the FDEs hold for $t \geq t_0$ almost everywhere with respect to the Lebesgue measure. The theorem was first proposed in our previous publication [23, Lemma 4], but the proofs were not provided in detail. Compared to the proofs of the conventional Krasovskiĭ stability theorem, the mean value theorem for vector-valued functions is no longer applicable for the FDEs with weak derivative, as the solution to FDEs is not differentiable everywhere. This problem is circumvented in this note by an application of the properties of Lebesgue integrals. The resulting stability theorem is almost identical to the conventional Krasovskiĭ stability theorem. However, it can be applied to FDEs that do not strictly satisfy the Marachkov boundedness condition [24, Theorem 6.1.3] and are not piecewise continuous in time. Such FDEs can be found in systems with time-varying delays [23], [25]–[27] that are integrable but not piecewise continuous.

*Notation*

Let $\mathbb{R}_+ := \{x \in \mathbb{R} : x > 0\}$ and $\mathbb{R}_{\geq a} := \{x \in \mathbb{R} : x \geq a\}$ with $a \in \mathbb{R}$ where $\mathbb{R}$ denotes the set of all real numbers. Standard $p$-norm for $\mathbb{R}^n$ is defined as $\mathbb{R}^n \ni \mathbf{x} \to \|\mathbf{x}\|_p := (\sum_{i=1}^n |x_i|^p)^{\frac{1}{p}}$ with $p \in \mathbb{N}$. $\mathcal{M}(\mathcal{X};\mathbb{R}^d)$ stands for the set containing all measurable functions defined from Lebesgue measurable set $\mathcal{X}$ to $\mathbb{R}^d$ endowed with the Borel algebra. We use $\mathcal{C}(\mathcal{X};\mathbb{R}^n)$ to denote the Banach space of continuous functions endowed with a uniform norm $\|\boldsymbol{f}(\cdot)\|_\infty := \sup_{\tau\in\mathcal{X}} \|\boldsymbol{f}(\tau)\|_2$, whereas $\mathcal{C}_\delta(\mathcal{X};\mathbb{R}^n) := \{\boldsymbol{f}(\cdot) \in \mathcal{C}(\mathcal{X};\mathbb{R}^n) : \|\boldsymbol{f}(\cdot)\|_\infty < \delta\}$ represents a normed bounded continuous functions space with $\delta > 0$. Lebesgue spaces are represented by $\mathcal{L}^p(\mathcal{X};\mathbb{R}^n) := \{\boldsymbol{f}(\cdot) \in \mathcal{M}(\mathcal{X};\mathbb{R}^n) : \|\boldsymbol{f}(\cdot)\|_p < +\infty\}$ with semi-norm $\|\boldsymbol{f}(\cdot)\|_p := \left(\int_{\mathcal{X}} \|\boldsymbol{f}(x)\|_2^p \mathsf{d}x\right)^{1/p}$. The space $\mathcal{K}_\infty$ of comparison functions follows the standard definition in [5, Chapter 4.4]. Notation $\widetilde{\forall} x \in \mathcal{X}, \mathtt{P}(x)$ means that property $\mathtt{P}(x)$ holds almost everywhere for $x \in \mathcal{X}$ w.r.t the Lebesgue measure. Symbol $\mathbf{0}_n$ represents an $n \times 1$ column vector and $\mathbb{0}_n(\cdot)$ corresponds to a "zero function" satisfying $\forall\theta \in [-r, 0]$, $\mathbb{0}_n(\theta) = \mathbf{0}_n$.

## II. Functional Differential Equation in the Extended Sense

Consider an FDE

$$\begin{aligned} &\widetilde{\forall} t \in \mathcal{T},\ \ \dot{\boldsymbol{x}}(t) = \boldsymbol{f}(t, \mathbf{x}_t(\cdot)),\ \ \emptyset \neq \mathcal{T} = [t_0, +\infty) \cap \mathcal{U}, \\ &\forall\theta \in [-r, 0],\ \ \boldsymbol{x}(t_0 + \theta) = \mathbf{x}_{t_0}(\theta) = \boldsymbol{\phi}(\theta),\ \ r > 0, \\ &\forall t \in \mathbb{R},\ \mathbf{0}_n = \boldsymbol{f}(t, \mathbb{0}_n(\cdot)),\ \forall\theta \in [-r, 0],\ \boldsymbol{x}(t + \theta) = \mathbf{x}_t(\theta) \end{aligned} \tag{1}$$

in the extended sense, where $t_0 \in \mathcal{U} \subseteq \mathbb{R}$ and the FDE holds for $t \in \mathcal{T}$ almost everywhere satisfying the initial condition $\boldsymbol{\phi}(\cdot) \in \mathcal{C}([-r, 0];\mathbb{R}^n)$. The existence of solutions to the FDE in (1) is guaranteed if $\boldsymbol{f} : \mathbb{R} \times \mathcal{X} \to \mathbb{R}^n$ satisfies the *Carathéodory condition* [13, section 2.6] on open set $\mathcal{X} \subseteq \mathcal{C}([-r, 0];\mathbb{R}^n)$. Namely, $\boldsymbol{f}(\cdot, \boldsymbol{\phi}(\cdot))$ is measurable for any given $\boldsymbol{\phi}(\cdot) \in \mathcal{X}$, and $\boldsymbol{f}(t, \bullet)$ is continuous for all $t \in \mathbb{R}$, and for any compact $\mathcal{Q} \subseteq \mathbb{R} \times \mathcal{X}$,

$$\exists m(\cdot) \in \mathcal{L}^1(\mathcal{Q}, \mathbb{R}),\ \forall(t, \boldsymbol{\phi}(\cdot)) \in \mathcal{Q},\ \|\boldsymbol{f}(t, \boldsymbol{\phi}(\cdot))\|_1 \leq m(t). \tag{2}$$

Finally, the uniqueness of the initial value problem in (1) is guaranteed if $\boldsymbol{f}(t, \boldsymbol{\phi}(\cdot))$ is integral locally Lipschitz in $\boldsymbol{\phi}(\cdot)$. Note that $\mathbf{x}_t(\cdot) \in \mathcal{C}([-r, 0];\mathbb{R}^n)$ in (1) is absolutely continuous for all $t > t_0 \in \mathbb{R}$, which means that $\boldsymbol{x}(\cdot)$ is differentiable almost everywhere on $[t_0, \infty]$.

The stabilities of FDEs [7], [13] can be defined similarly to the stabilities of ODEs [3], [5]. A diagram illustrating the notions of stability is presented in Figure 1.

*Definition 1:* Assume $\boldsymbol{x}(t)$ is the unique solution to the initial value problem in (1). The origin of the FDE in (1) is

***Lyapunov Stable*** in region $\mathcal{Y} \subseteq \mathcal{C}([-r, 0];\mathbb{R}^n)$, if $\forall\epsilon > 0$, $\forall t_0 \in \mathbb{R}$, $\exists\delta(t_0, \epsilon) > 0$, $\forall\boldsymbol{\phi}(\cdot) \in \mathcal{C}_{\delta(t_0,\epsilon)}([-r, 0];\mathbb{R}^n) \cap \mathcal{Y}$, $\forall t \geq t_0$, $\|\mathbf{x}_t(\cdot)\|_\infty < \epsilon$.

***Uniform Stable*** in region $\mathcal{Y} \subseteq \mathcal{C}([-r, 0];\mathbb{R}^n)$, if $\forall\epsilon > 0$, $\exists\delta(\epsilon) > 0$, $\forall t_0 \in \mathbb{R}$, $\forall\boldsymbol{\phi}(\cdot) \in \mathcal{C}_{\delta(\epsilon)}([-r, 0];\mathbb{R}^n) \cap \mathcal{Y}$, $\forall t \geq t_0$, $\|\mathbf{x}_t(\cdot)\|_\infty < \epsilon$.

**Uniformly Asymptotically stable** in $\mathcal{Y} \subseteq \mathcal{C}([-r, 0];\mathbb{R}^n)$, if it is uniformly stable in region $\mathcal{Y}$, and $\exists\delta > 0$, $\forall\eta > 0$, $\exists\beta(\eta) \geq 0$, $\forall\boldsymbol{\phi}(\cdot) \in \mathcal{C}_\delta([-r, 0];\mathbb{R}^n) \cap \mathcal{Y}$, $\forall t_0 \in \mathbb{R}$, $\forall t \geq t_0 + \beta(\eta)$, $\|\mathbf{x}_t(\cdot)\|_\infty < \eta$.

**Globally Uniformly Asymptotically stable** for any initial value $\boldsymbol{\phi}(\cdot) \in \mathcal{C}([-r, 0];\mathbb{R}^n)$, if it is uniformly stable and $\forall\eta > 0$, $\forall\delta > 0$, $\exists\beta(\eta, \delta) > 0$, $\forall\boldsymbol{\phi}(\cdot) \in \mathcal{C}_\delta([-r, 0];\mathbb{R}^n)$, $\forall t_0 \in \mathbb{R}$, $\forall t \geq t_0 + \beta(\eta, \delta)$, $\|\mathbf{x}_t(\cdot)\|_\infty < \eta$.

Given the definition of both FDEs in the extended sense and the stabilities of their equilibrium points, the main theorem in this note is enunciated as follows.

*Theorem 1:* Consider the FDE in (1) that satisfies both the *Carathéodory condition* and integral locally Lipschitz condition for ensuring the existence and uniqueness of the solutions. Moreover, we assume $\forall\delta > 0$, $\exists c > 0$,

$$\forall\boldsymbol{\phi}(\cdot) \in \mathcal{C}_\delta([-r, 0];\mathbb{R}^n),\ \widetilde{\forall} t \in \mathbb{R},\ \|\boldsymbol{f}(t, \boldsymbol{\phi}(\cdot))\|_1 < c. \tag{3}$$

Let $\mathcal{U} = \mathbb{R}$. Then the trivial solution $\boldsymbol{x}(t) \equiv \mathbf{0}_n$ of (1) is globally uniformly asymptotically stable if there exist three functions $\alpha_1(\cdot); \alpha_2(\cdot); \alpha_3(\cdot) \in \mathcal{K}_\infty$, and a continuous functional $\mathsf{v} : \mathbb{R} \times \mathcal{C}([-r, 0];\mathbb{R}^n) \to \mathbb{R}$ such that

$$\forall t \in \mathbb{R},\ \alpha_1(\|\boldsymbol{\phi}(0)\|_2) \leq \mathsf{v}(t, \boldsymbol{\phi}(\cdot)) \leq \alpha_2(\|\boldsymbol{\phi}(\cdot)\|_\infty) \tag{4}$$

$$\widetilde{\forall} t \geq t_0 \in \mathbb{R}, \qquad \frac{\mathsf{d}}{\mathsf{d}t}\mathsf{v}(t, \mathbf{x}_t(\cdot)) \leq -\alpha_3(\|\boldsymbol{x}(t)\|_2) \tag{5}$$

for all initial condition $\boldsymbol{\phi}(\cdot) \in \mathcal{C}([-r, 0];\mathbb{R}^n)$ in (1), where $\|\boldsymbol{\phi}(\cdot)\|_\infty := \max_{-r_2\leq\tau\leq 0} \|\boldsymbol{\phi}(\tau)\|$, and $\mathbf{x}_t(\cdot)$, $\boldsymbol{x}(\cdot)$ satisfy the FDE in (5) for almost all $t \geq t_0$. Moreover, $\mathcal{K}_\infty$ functions follow the standard definition in [5]. Notation $\widetilde{\forall}$ denotes for almost all with respect to the Lebesgue measure [17].

*Proof:* The proof here is based on the procedure in [13, Chapter 5, Theorem 2.1] with additional steps and modifications that address the extended sense of FDEs in (1). To establish the uniform stability of the origin, let

$$\mathbb{R}_{\geq 0} \ni \epsilon \mapsto \delta(\epsilon) = 1/2 \min\left(\epsilon, \alpha_2^{-1}(\alpha_1(\epsilon))\right) \tag{6}$$

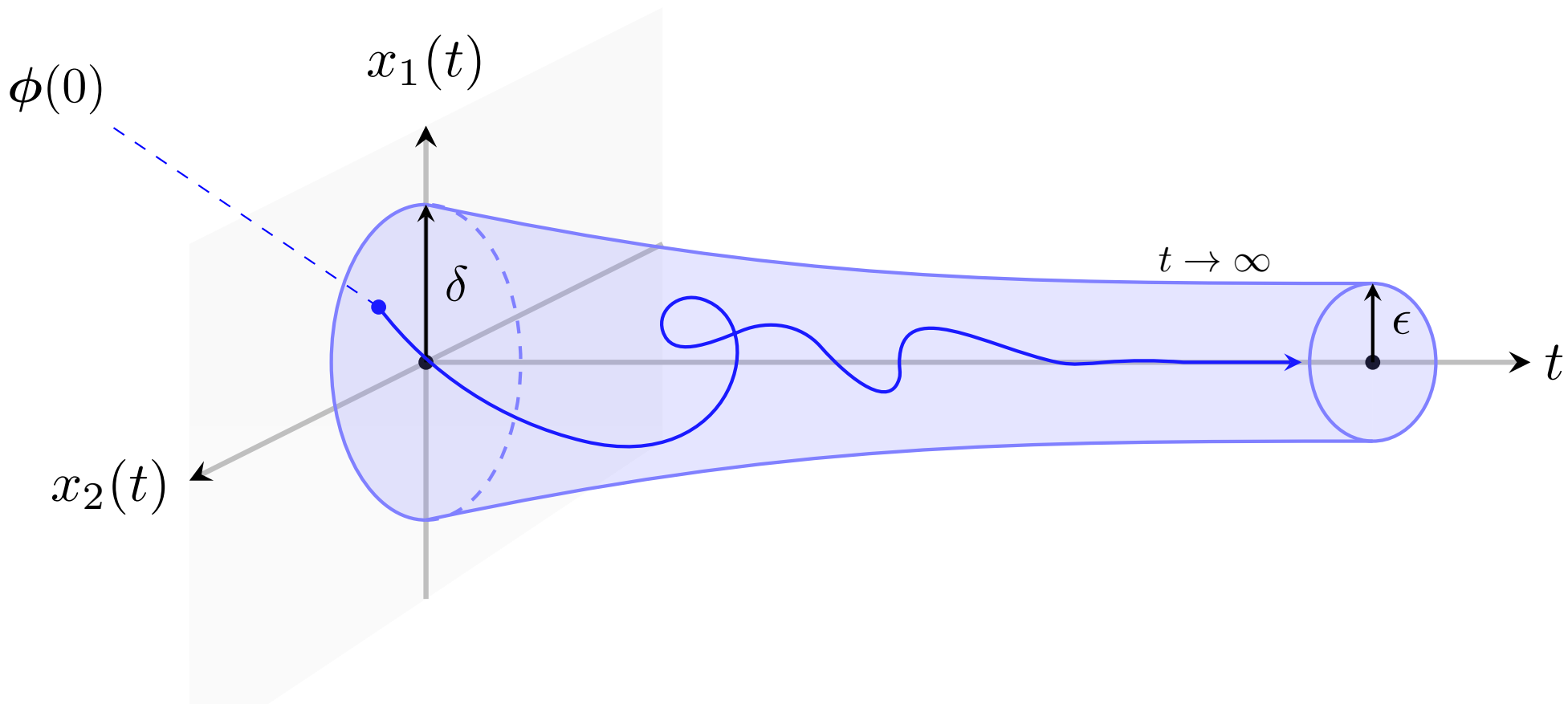

Fig. 1: Illustration of stability concepts using the $\epsilon - \delta$ paradigm, created with code from `https://tex.stackexchange.com/a/560965/104839`

where $\alpha_2^{-1}(\cdot)$ is well defined since $\alpha_2(\cdot) \in \mathcal{K}_\infty$. It is evident that $\delta(\cdot) \in \mathcal{K}_\infty$ and satisfies $\forall \epsilon > 0,\ 0 < \delta(\epsilon) < \epsilon$ and $\delta(\epsilon) < \alpha_2^{-1}(\alpha_1(\epsilon))$ that further implies

$$\forall \epsilon > 0, \ \ \alpha_2(\delta(\epsilon)) < \alpha_1(\epsilon) \tag{7}$$

since $\alpha_2(\cdot) \in \mathcal{K}_\infty$. By (5), it holds true that $\widetilde{\forall} t \geq t_0 \in \mathbb{R},\ \frac{\mathsf{d}}{\mathsf{d}t}\mathsf{v}(t, \mathbf{x}_t(\cdot)) \leq 0$. Now, by applying the Fundamental Theorem of Calculus for Lebesgue integrals [17] to the previous proposition, we see that

$$\begin{aligned} &\forall t \geq t_0, \forall \boldsymbol{\phi}(\cdot) \in \mathcal{C}([-r,0];\mathbb{R}^n),\ \int_{t_0}^{t} \frac{\mathsf{d}}{\mathsf{d}\tau}\mathsf{v}(\tau, \mathbf{x}_\tau(\cdot))\mathsf{d}\tau \\ &= \mathsf{v}(t, \mathbf{x}_t(\cdot)) - \mathsf{v}(t_0, \mathbf{x}_{t_0}(\cdot)) = \mathsf{v}(t, \mathbf{x}_t(\cdot)) - \mathsf{v}(t_0, \boldsymbol{\phi}(\cdot)) \leq 0 \end{aligned} \tag{8}$$

which further implies that $\forall t_0 \in \mathbb{R}, \forall t \geq t_0, \forall \epsilon > 0, \forall \boldsymbol{\phi}(\cdot) \in \mathcal{C}_{\delta(\epsilon)}([-r,0];\mathbb{R}^n)$, we have

$$\begin{aligned} \alpha_1(\|\boldsymbol{x}(t)\|_2) \leq \mathsf{v}(t, \mathbf{x}_t(\cdot)) \leq \mathsf{v}(t_0, \boldsymbol{\phi}(\cdot)) \leq \alpha_2(\|\boldsymbol{\phi}(\cdot)\|_\infty) \\ < \alpha_2(\delta(\epsilon)) < \alpha_1(\epsilon) \end{aligned} \tag{9}$$

by the relations in (4) and (7), where $\mathcal{C}_{\delta(\epsilon)}([-r,0];\mathbb{R}^n) := \{\boldsymbol{\phi}(\cdot) \in \mathcal{C}([-r,0];\mathbb{R}^n) : \|\boldsymbol{\phi}(\cdot)\|_\infty < \delta(\epsilon)\}$. Therefore, for all $\epsilon > 0$, and $\boldsymbol{\phi}(\cdot) \in \mathcal{C}_{\delta(\epsilon)}([-r,0];\mathbb{R}^n)$, we have

$$\forall t_0 \in \mathbb{R},\ \forall t \geq t_0,\ \|\boldsymbol{x}(t)\|_2 < \epsilon \quad \text{since} \ \ \alpha_1(\cdot) \in \mathcal{K}_\infty, \tag{10}$$

where $\delta(\epsilon) = 1/2 \min\left(\epsilon, \alpha_2^{-1}(\alpha_1(\epsilon))\right)$ is independent of $t_0 \in \mathbb{R}$ and $\lim_{\epsilon \to +\infty} \delta(\epsilon) = +\infty$ since $\delta(\cdot) \in \mathcal{K}_\infty$. Now one can further infer $\forall \epsilon > 0,\ \exists \delta > 0,\ \forall \boldsymbol{\phi}(\cdot) \in \mathcal{C}_\delta([-r,0];\mathbb{R}^n)$

$$\forall t_0 \in \mathbb{R},\ \forall t \geq t_0,\ \|\mathbf{x}_t(\cdot)\|_\infty \leq \max_{\tau \geq t_0} \|\boldsymbol{x}(\tau)\| < \epsilon \tag{11}$$

form (10) which demonstrates uniform stability.

To prove global uniform asymptotic stability, we employ proof by contradiction. Note that the origin is globally uniform asymptotically stable if it is uniformly stable and $\forall \eta > 0, \forall \delta > 0,\ \exists \beta(\eta, \delta) \geq 0,\ \forall \boldsymbol{\phi}(\cdot) \in \mathcal{C}_\delta([-r,0];\mathbb{R}^n)$

$$\forall t_0 \in \mathbb{R},\ \forall t \geq t_0 + \beta(\eta, \delta),\ \|\mathbf{x}_t(\cdot)\|_\infty < \eta. \tag{12}$$

Now let us assume that the origin is uniformly stable but there exist $\epsilon > 0$ and $\delta > 0$ and function $\boldsymbol{\phi}(\cdot) \in \mathcal{C}_\delta([-r,0];\mathbb{R}^n)$ and $t_0 \in \mathbb{R}$ such that

$$\forall t \geq t_0,\ \|\mathbf{x}_t(\cdot)\|_\infty \geq \epsilon. \tag{13}$$

Given the definition $\|\mathbf{x}_t(\cdot)\|_\infty := \max_{\tau \in [-r,0]} \|\boldsymbol{x}(t+\tau)\|_2$ with (13), it implies there exist constants $\epsilon > 0,\ \delta > 0$ and function $\boldsymbol{\phi}(\cdot) \in \mathcal{C}_\delta([-r,0];\mathbb{R}^n)$ and $t_0 \in \mathbb{R}$ such that

$$\forall t \geq t_0,\ \exists \lambda \in [t-r, t],\ \|\boldsymbol{x}(\lambda)\|_2 \geq \epsilon. \tag{14}$$

Let $\epsilon > 0,\ \delta > 0,\ \boldsymbol{\phi}(\cdot) \in \mathcal{C}_\delta([-r,0];\mathbb{R}^n)$ and $t_0 \in \mathbb{R}$ such that (14) is true, then we see that there exists a sequence $\mathbb{N} \ni k \to t_k \in [t_0, \infty)$ such that

$$\forall k \in \mathbb{N},\ (2k-1)r \leq t_k - t_0 \leq 2kr \ \&\ \|\boldsymbol{x}(t_k)\|_2 \geq \epsilon. \tag{15}$$

Additionally,

$$\begin{aligned} \|\boldsymbol{x}(t)\|_2 &= \left\| \boldsymbol{x}(t_k) + \int_{t_k}^{t} \dot{\boldsymbol{x}}(\tau)\mathsf{d}\tau \right\|_2 \\ &\underbrace{\geq \|\boldsymbol{x}(t_k)\|_2 - \left\| \int_{t_k}^{t} \dot{\boldsymbol{x}}(\tau)\mathsf{d}\tau \right\|_2}_{\text{reverse triangle inequality}} \\ &= \|\boldsymbol{x}(t_k)\|_2 - \left\| \int_{t_k}^{t} \boldsymbol{f}(\tau, \mathbf{x}_\tau(\cdot))\mathsf{d}\tau \right\|_2 \\ &\underbrace{\geq \|\boldsymbol{x}(t_k)\|_2 - \left\| \int_{t_k}^{t} \boldsymbol{f}(\tau, \mathbf{x}_\tau(\cdot))\mathsf{d}\tau \right\|_1}_{\|\cdot\|_1 \geq \|\cdot\|_2} \\ &= \|\boldsymbol{x}(t_k)\|_2 - \sum_{i=1}^{n} \left| \int_{t_k}^{t} f_i(\tau, \mathbf{x}_\tau(\cdot))\mathsf{d}\tau \right| \\ &\underbrace{\geq \|\boldsymbol{x}(t_k)\|_2 - \left| \sum_{i=1}^{n} \int_{t_k}^{t} |f_i(\tau, \mathbf{x}_\tau(\cdot))|\,\mathsf{d}\tau \right|}_{\sum_i |\int f_i(\tau)\mathsf{d}\tau| \leq \sum_i \int |f_i(\tau)|\mathsf{d}\tau} \\ &= \|\boldsymbol{x}(t_k)\|_2 - \left| \int_{t_k}^{t} \sum_{i=1}^{n} |f_i(\tau, \mathbf{x}_\tau(\cdot))|\,\mathsf{d}\tau \right| \\ &= \|\boldsymbol{x}(t_k)\|_2 - \left| \int_{t_k}^{t} \|\boldsymbol{f}(\tau, \mathbf{x}_\tau(\cdot))\|_1\,\mathsf{d}\tau \right| \end{aligned} \tag{16}$$

holds true for all $t \geq t_0$ and $k \in \mathbb{N}$ based on the properties of Lebesgue integrals and norms. Since the origin is uniformly stable as in (11) and $\delta(\cdot)$ in (6) is surjective, we see that $\forall\delta>0,\ \exists\varepsilon(\delta)>0$ such that $\forall\boldsymbol{\phi}(\cdot)\in\mathcal{C}_\delta([-r,0];\mathbb{R}^n)$

$$\forall t_0\in\mathbb{R},\ \forall t\geq t_0,\ \|\mathbf{x}_t(\cdot)\|_\infty<\varepsilon(\delta)$$

by which we know that $\forall\delta>0,\ \exists c(\varepsilon(\delta))>0$,

$$\forall\boldsymbol{\phi}(\cdot)\in\mathcal{C}_\delta\left([-r,0];\mathbb{R}^n\right),\ \widetilde{\forall} t\in\mathbb{R},\ \|\boldsymbol{f}\left(t,\mathbf{x}_t(\cdot)\right)\|_1<c \quad (17)$$

following from (3). Now consider proposition $\forall t\geq t_0,\ \forall k\in\mathbb{N},\ \left|\int_{t_k}^t\|\boldsymbol{f}(\tau,\mathbf{x}_\tau(\cdot))\|_1\,\mathsf{d}\tau\right|<\left|\int_{t_k}^t c\,\mathsf{d}\tau\right|=c|t-t_k|$ as a result of (17) under given $\delta>0$ and $\boldsymbol{\phi}(\cdot)\in\mathcal{C}_\delta([-r,0];\mathbb{R}^n)$, we can then conclude

$$\begin{aligned}\|\boldsymbol{x}(t)\|_2&\geq\|\boldsymbol{x}(t_k)\|_2-\left|\int_{t_k}^t\|\boldsymbol{f}(\tau,\mathbf{x}_\tau(\cdot))\|_1\,\mathsf{d}\tau\right|\\&>\|\boldsymbol{x}(t_k)\|_2-\left|\int_{t_k}^t c\,\mathsf{d}\tau\right|=\|\boldsymbol{x}(t_k\|_2-c|t-t_k|\\&\geq\epsilon-c\frac{\epsilon}{2c}=\frac{\epsilon}{2}\end{aligned}\quad (18)$$

for all $k\in\mathbb{N}$ and $t\in\mathcal{T}_k:=\left[t_k-\frac{\epsilon}{2c},t_k+\frac{\epsilon}{2c}\right]$. Consequently, we have

$$\begin{aligned}&\widetilde{\forall} t\in\left[\mathbb{R}_{\geq t_0}\cap\bigcup_{k\in\mathbb{N}}\mathcal{T}_k\right],\ \frac{\mathsf{d}}{\mathsf{d}t}\mathsf{v}(t,\mathbf{x}_t(\cdot))\leq-\alpha_3\left(\epsilon/2\right)\\&\widetilde{\forall} t\in\mathbb{R}_{\geq t_0},\ \ \frac{\mathsf{d}}{\mathsf{d}t}\mathsf{v}(t,\mathbf{x}_t(\cdot))\leq 0.\end{aligned}\quad (19)$$

Since $c(\varepsilon(\delta))>0$ in $\mathcal{T}_k=[t_k-\epsilon/2c,t_k+\epsilon/2c]$ can be made arbitrarily large for any $\delta>0$, we can assume that intervals $\mathcal{T}_k$ are pairwise disjoint (no overlaps) $\bigcup_{i,j\in\mathbb{N},i\neq j}\mathcal{T}_i\cap\mathcal{T}_j=\emptyset$ and $t_1-\epsilon/2c\geq t_0$. As a result, we have

$$\begin{aligned}&\mathsf{v}(t_k,\mathbf{x}_{t_k}(\cdot))-\mathsf{v}(t_0,\boldsymbol{\phi}(\cdot))=\int_{t_0}^{t_k}\frac{\mathsf{d}}{\mathsf{d}\tau}\mathsf{v}(\tau,\mathbf{x}_\tau(\cdot))\mathsf{d}\tau\\&=\int_{\bigcup_{i=1}^{k-1}\mathcal{T}_i}\frac{\mathsf{d}}{\mathsf{d}\tau}\mathsf{v}(\tau,\mathbf{x}_\tau(\cdot))\mathsf{d}\tau+\underbrace{\int_{[t_0,t_k]\setminus\bigcup_{i=1}^{k-1}\mathcal{T}_i}\frac{\mathsf{d}}{\mathsf{d}\tau}\mathsf{v}(\tau,\mathbf{x}_\tau(\cdot))\mathsf{d}\tau}_{\leq 0}\\&\leq-\int_{\bigcup_{i=1}^{k-1}\mathcal{T}_i}\alpha_3\left(\epsilon/2\right)\mathsf{d}\tau=-\sum_{i=1}^{k-1}\int_{\mathcal{T}_i}\alpha_3\left(\epsilon/2\right)\mathsf{d}\tau\\&=-\alpha_3\left(\epsilon/2\right)\frac{\epsilon}{c}(k-1),\forall k\in\mathbb{N}\end{aligned}\quad (20)$$

by (19). This further indicates that for all $k\in\mathbb{N}$ we have

$$\begin{aligned}\mathsf{v}(t_k,\mathbf{x}_{t_k}(\cdot))&\leq\mathsf{v}(t_0,\boldsymbol{\phi}(\cdot))-\alpha_3\left(\epsilon/2\right)\frac{\epsilon}{c}(k-1)\\&\leq\alpha_2\left(\|\boldsymbol{\phi}(\cdot)\|_\infty\right)-\alpha_3\left(\epsilon/2\right)\frac{\epsilon}{c}(k-1)\\&<\alpha_2(\delta)-\alpha_3\left(\epsilon/2\right)\frac{\epsilon}{c}(k-1)\end{aligned}\quad (21)$$

by (4) and the fact that $\|\boldsymbol{\phi}(\cdot)\|_\infty<\delta$ and $\alpha_2(\cdot)\in\mathcal{K}_\infty$. Now it is easy to see that

$$\alpha_2(\delta)-\alpha_3\Big(\frac{\epsilon}{2}\Big)\frac{\epsilon}{c}(k-1)<0\iff\frac{\alpha_2(\delta)c}{\alpha_3\left(\epsilon/2\right)\epsilon}+1<k.$$

Let $\kappa(\epsilon,\delta)=\left\lceil\frac{\alpha_2(\delta)c(\varepsilon(\delta))}{\alpha_3(\epsilon/2)\epsilon}\right\rceil+1$, then we can obtain $\forall k>\kappa(\epsilon,\delta)$, $\mathsf{v}(t_k,\mathbf{x}_{t_k}(\cdot))<0$ by (21), which is a contradiction of (4). In consequence, (15) cannot be true for $t_k$ with any $k>\kappa(\epsilon,\delta)$, which means that $\exists k\leq\kappa(\epsilon,\delta),\ \|\mathbf{x}_{t_k}(\cdot)\|_\infty<\epsilon$. This further indicates that for all $\epsilon>0,\ \delta>0$ and $\boldsymbol{\phi}(\cdot)\in\mathcal{C}_\delta([-r,0];\mathbb{R}^n)$, we have

$$\forall t_0\in\mathbb{R},\ \exists\theta\in[t_0,t_0+2\kappa(\epsilon,\delta)r],\ \|\mathbf{x}_\theta(\cdot)\|_\infty<\epsilon \quad (22)$$

given the proposition in (15) and the fact that

$$[t_0+(2\kappa(\epsilon,\delta)-1)r,t_0+2\kappa(\epsilon,\delta)r]\subset[t_0,t_0+2\kappa(\epsilon,\delta)r].$$

Now let $\epsilon>0$ in (22) be

$$\epsilon(\eta)=1/3\min\left(\eta,\alpha_2^{-1}\left(\alpha_1(\eta)\right)\right) \quad (23)$$

with a given $\eta>0$, and assume $\boldsymbol{\phi}(\cdot),\ t_0,\ \theta$ in (22) are also given. Note that the structure of $\epsilon(\cdot)$ in (23) guarantees $\epsilon(\cdot)\in\mathcal{K}_\infty$ and $\alpha_2(\epsilon(\eta))<\alpha_1(\eta)$ for any $\eta>0$ similar to the property in (7). From the computation in (8) and (5), we also know that for all $t\geq\theta$ and $\boldsymbol{\phi}(\cdot)\in\mathcal{C}([-r,0];\mathbb{R}^n)$

$$\int_\theta^t\frac{\mathsf{d}}{\mathsf{d}\tau}\mathsf{v}(\tau,\mathbf{x}_\tau(\cdot))=\mathsf{v}\left(t,\mathbf{x}_t(\cdot)\right)-\mathsf{v}\big(\theta,\mathbf{x}_\theta(\cdot)\big)\leq 0. \quad (24)$$

With (4) and (22)–(24) and $\alpha_2(\epsilon(\eta))<\alpha_1(\eta)$, we find that $\forall\eta>0,\ \forall\delta>0,\ \forall t_0\in\mathbb{R},\ \forall\boldsymbol{\phi}(\cdot)\in\mathcal{C}_\delta([-r,0];\mathbb{R}^n)$,

$$\begin{aligned}\forall t\geq\theta,\ \ \alpha_1(\|\boldsymbol{x}(t)\|_2)\leq\mathsf{v}(t,\mathbf{x}_t(\cdot))&\leq\mathsf{v}\left(\theta,\mathbf{x}_\theta(\cdot)\right)\\&\leq\alpha_2(\|\mathbf{x}_\theta(\cdot)\|_\infty)\\&<\alpha_2(\epsilon(\eta))\\&<\alpha_1(\eta).\end{aligned}\quad (25)$$

Since $\theta\leq t_0+2\kappa(\epsilon,\delta)r$ in (22), relation (25) further implies

$$\forall t\geq t_0+2\kappa(\epsilon,\delta)r\geq\theta,\quad\alpha_1(\|\boldsymbol{x}(t)\|_2)<\alpha_1(\eta)$$

and then $\|\boldsymbol{x}(t)\|_2<\eta$ since $\alpha_1(\cdot)\in\mathcal{K}_\infty$. Given that $2r\kappa(\epsilon(\eta),\delta)$ is independent of $t_0$, it proves the global asymptotic stability in (12) with $\beta=2\kappa(\epsilon,\delta)r$. This proves the global uniform asymptotic stability of $\boldsymbol{x}(t)\equiv\mathbf{0}_n$, as the uniform stability has been proved with $\delta(\cdot)\in\mathcal{K}_\infty$ in (6) satisfying $\lim_{\epsilon\to+\infty}\delta(\epsilon)=+\infty$. ■

We have purposely assumed $\mathsf{v}(\cdot,\bullet)$ is continuous in both arguments before the condition in (5) to ensure that no additional conservatism is imposed on $\mathsf{v}(\cdot,\bullet)$. Note that $\mathsf{v}(\cdot,\bullet)$ being continuous does not guarantee that $\dot{\mathsf{v}}(t,\mathbf{x}_t(\cdot))$ exists for almost all $t\geq t_0$. Since (5) is part of the statement in Theorem 1 to be ascertained, it is unnecessary to specify the conditions required to ensure $\dot{\mathsf{v}}(t,\mathbf{x}_t(\cdot))$ exists for almost all $t\geq t_0$. Meanwhile, it is vitally important to stress that being locally absolutely continuous in $t$ and locally Lipschitz in the second argument cannot ensure the derivative in (5) exists for $t\geq t_0$ almost everywhere, even though $\boldsymbol{x}(t)$ is locally absolutely continuous for all $t\geq t_0$ according to the structures in (1). This is because $\mathbf{x}_t(\cdot)$ always overlaps with some segment of the initial condition $\boldsymbol{\phi}(\cdot)\in\mathcal{C}\left([-r,0];\mathbb{R}^n\right)$, which can be totally non-differentiable, when $t_0\leq t\leq t_0+r$. For instance, if $\mathsf{v}(t,\mathbf{x}_t(\cdot))$ includes quadratic term $\boldsymbol{x}^\top(t-r)Q\boldsymbol{x}(t-r)$, then $\dot{\mathsf{v}}(t,\mathbf{x}_t(\cdot))$ cannot exist on $t\in[t_0,t_0+r]$ unless $\boldsymbol{\phi}(\cdot)\in\mathcal{W}^{1,2}\left([-r,0];\mathbb{R}^n\right)\subset\mathcal{C}\left([-r,0];\mathbb{R}^n\right)$. This noteworthy property was first studied by Pierdomenico Pepe in [28]–[30] from the perspective of input-to-state stability.

*Remark 1:* The condition in (3) is added to prove Theorem 1, which indicates that $\boldsymbol{f}(t,\boldsymbol{\phi}(\cdot))$ must be locally bounded in $\boldsymbol{\phi}(\cdot)$ and essentially bounded in $t$. This condition can be viewed as a measure-theoretic counterpart of the Marachkov boundedness condition [24, Theorem 6.1.3] required in the proof of the conventional Krasovskiĭ stability theorem for FDEs with standard derivatives.

*Remark 2:* Theorem 1 addresses the global asymptotic stability of (1). Non-global versions of various stability concepts of FDEs can be established following the proof of Theorem 1.

The conventional negative condition for the existence of Krasovskiĭ functional $\mathsf{v}(\cdot, \bullet)$ in [13, Chapter 5, Theorem 2.1] is defined using the upper right-hand Dini derivative [31]

$$\dot{\mathsf{v}}(t, \boldsymbol{\phi}(\cdot)) = \limsup_{h \to 0^+} \frac{\mathsf{v}(t+h, \mathbf{x}_{t+h}(\cdot)) - \mathsf{v}(t, \boldsymbol{\phi}(\cdot))}{h}.$$

To ensure $\dot{\mathsf{v}}(t, \boldsymbol{\phi}(\cdot))$ is finite for all $t \geq t_0$, which is particularly meaningful in the control contexts, functional $\mathsf{v}(\cdot, \bullet)$ must be locally Lipschitz in both arguments [13, Chapter 5, Theorem 2.1], which implies that $\mathsf{v}(t, \mathbf{x}_t(\cdot))$ is differentiable almost everywhere for $t > t_0$. (See Rademacher's Theorem [17]) In contrast, the negative condition in (5) is formulated using the language of measure theory instead of Dini derivatives. In the conventional formulation of FDEs using standard derivatives, $\boldsymbol{f}(t, \mathbf{x}_t(\cdot))$ in [13, Chapter 5, Theorem 2.1] must be continuous in $t$, hence $\boldsymbol{x}(t)$ has to be (left) differentiable. This property is utilized in the proof of the traditional stability theorem [13, Chapter 5, Theorem 2.1], [18, Theorem 1.3] where the mean value theorem is applied to $\|\boldsymbol{x}(t)\|$ to construct the corresponding step in (18). However, we can no longer apply this strategy in our proof since $\boldsymbol{x}(t)$ is not strictly differentiable for all $t \geq t_0$. This demonstrates the necessity of establishing the stability criteria in Theorem 1 as a separate stability theorem.

*Remark 3:* If we consider the equation in (1) as a switching system in $t$, then its stability can be analyzed using [32, Theorem 11] with a measurable switching signal $\sigma(t) = t$. A critical component of the proof of [32, Theorem 11] is the conclusion on input-to-state stability in [32, Theorem 9] that shows the equivalence between M-ISS and PC-ISS.

## III. Application to Systems with Time-Varying Delays and Disturbances

If $\boldsymbol{f}(\cdot, \bullet)$ in (1) is continuous or piecewise continuous in $t$ and locally Lipschitz in the second argument, then Theorem 1 becomes identical to the conventional Krasovskiĭ stability theorem [13, Section 5.2]. However, such continuity is not always guaranteed for systems operating in a real environment. For example, some networked control problems can be described by a linear system

$$\dot{\boldsymbol{x}}(t) = A_1 \boldsymbol{x}(t) + A_2 \boldsymbol{x}(t - r(t)), \ \boldsymbol{\phi}(\cdot) \in \mathcal{C}([-r_2, 0]; \mathbb{R}) \quad (26)$$

with $A_1, A_2 \in \mathbb{R}^{n \times n}$ and time-varying delay $r(t)$ that is **neither differentiable everywhere nor piecewise continuous.** In many cases, $r(\cdot) \in \mathcal{L}^\infty(\mathbb{R}; \mathbb{R})$ is integrable and essentially bounded satisfying $\widetilde{\forall} t \in \mathbb{R}$, $r_1 \leq r(t) \leq r_2$, which can exhibit a countablely infinite number of jumps with arbitrary values. For example, let $r(t) = \eta(t) + n(t)\mathbb{1}_{\mathbb{Q}}(t)$ where $\eta(t)$ is any piecewise continuous function in $t$ and $n(t)$ is a continuous time stochastic process with indicator function $\mathbb{1}_{\mathbb{Q}}(t)$ over rational number set $\mathbb{Q}$. Since $r(t)$ may not be continuous but only integrable, the analysis of (26) cannot be achieved using the traditional Krasovskiĭ stability theorem, even if we were to consider the novel Krasovskiĭ stability theorems [14] proposed for time-varying delay systems. This is because the right hand side of $\dot{\boldsymbol{x}}(t)$ in (26) with $r(\cdot) \in \mathcal{L}^\infty(\mathbb{R}; \mathbb{R})$ is neither continuous nor piecewise continuous as required by the prerequisites of the Krasovskiĭ stability theorems in [13, Section 5.2] or [14], [15]. On the other hand, if we interpret (26) as

$$\widetilde{\forall} t \geq t_0, \ \dot{\boldsymbol{x}}(t) = A_1 \boldsymbol{x}(t) + A_2 \boldsymbol{x}(t - r(t)) \quad (27)$$

where the FDE holds for $t \geq t_0$ almost everywhere, then the extended differential equation with $r(\cdot) \in \mathcal{L}^\infty(\mathbb{R}; \mathbb{R})$ is well-defined.

Utilizing the properties of Lebesgue-Stieltjes integral, (27) can be rewritten as

$$\widetilde{\forall} t \geq t_0, \dot{\boldsymbol{x}}(t) = A_1 \boldsymbol{x}(t) + \int_{-r_2}^{0} A_2 \boldsymbol{x}(t+\tau) \mathsf{d} \mathit{1}(\tau + r(t)) \quad (28)$$

where $\mathit{1}(\cdot)$ is the standard Heaviside step function. It is clear to see that the right-hand side of (28) is a time-varying linear operator that satisfies the *Carathéodory conditions* and be locally Lipschitz, and also satisfies the boundedness condition in (3). Consequently, we can apply Theorem 1 to analyze the stability of the origin of (28), exactly like we can employ the conventional Krasovskiĭ stability theorem to (28) if $r(\cdot) \in \mathcal{C}(\mathbb{R}; [r_1, r_2])$. This shows that the existing methods [26], [33], [34] for the stability analysis of (28) remain valid, as the conclusions are supported by Theorem 1 without a need to change the derivation procedures.

Another representative example of an application of Theorem 1 can be found in [23], where a distributed delay system

$$\begin{aligned}
\dot{\boldsymbol{x}}(t) &= A_1 \boldsymbol{x}(t) + \int_{-r(t)}^{0} \widetilde{A}_2(\tau) \boldsymbol{x}(t+\tau) \mathsf{d}\tau \\
&\quad + B_1 \boldsymbol{u}(t) + \int_{-r(t)}^{0} \widetilde{B}_2(\tau) \boldsymbol{u}(t+\tau) \mathsf{d}\tau + D_1 \boldsymbol{w}(t), \\
\boldsymbol{z}(t) &= C_1 \boldsymbol{x}(t) + \int_{-r(t)}^{0} \widetilde{C}_2(\tau) \boldsymbol{x}(t+\tau) \mathsf{d}\tau \\
&\quad + B_4 \boldsymbol{u}(t) + \int_{-r(t)}^{0} \widetilde{B}_5(\tau) \boldsymbol{u}(t+\tau) \mathsf{d}\tau + D_2 \boldsymbol{w}(t), \\
&\forall \theta \in [-r_2, 0], \ \ \boldsymbol{x}(t_0 + \theta) = \boldsymbol{\phi}(\theta)
\end{aligned} \quad (29)$$

is investigated with time-varying delay $r(\cdot)$ that satisfies the same constraints as those in (26), and $\boldsymbol{w}(\cdot) \in \mathcal{L}^2([t_0, \infty); \mathbb{R}^n)$. Clearly, the system in (29) satisfies the same prerequisites for applying Theorem 1, as the operators in (26) are all linear. Moreover, the presence of $\boldsymbol{w}(\cdot) \in \mathcal{L}^2([t_0, \infty); \mathbb{R}^n)$ implies that the FDE in (29) cannot hold for all $t \geq t_0$ even if $r(\cdot)$ is strictly continuous. Thus an FDE with non-continuous disturbance $\boldsymbol{w}(\cdot)$ must be interpreted in the extended sense. This is actually another significant advantage of the Carathéodory framework that stability analysis and performance objectives with disturbances such as dissipativity can be addressed simultaneously using a single framework. Specifically, dissipativity can be defined as:

*Definition 2:* The time delay system in (29) with supply rate function $\mathsf{s}(\boldsymbol{z}(t), \boldsymbol{w}(t))$ is dissipative if there exists a continuous functional $\mathsf{v}(\cdot, \bullet) : \mathbb{R} \times \mathcal{C}([-r_2, 0]; \mathbb{R}^n) \to \mathbb{R}$ such that

$$\forall t \geq t_0, \ \ \mathsf{v}(t, \mathbf{x}_t(\cdot)) - \mathsf{v}(t_0, \mathbf{x}_{t_0}(\cdot)) \leq \int_{t_0}^{t} \mathsf{s}(\boldsymbol{z}(\theta), \boldsymbol{w}(\theta)) \mathsf{d}\theta \quad (30)$$

where $t_0 \in \mathbb{R}$ and regulate output $\boldsymbol{z}(t)$ and disturbance $\boldsymbol{w}(t)$ are given by (29). Moreover, $\mathbf{x}_t(\cdot)$ in (31) is defined by the proposition $\forall t \geq t_0$, $\forall \theta \in [-r_2, 0]$, $\mathbf{x}_t(\theta) = \boldsymbol{x}(t + \theta)$ where $\boldsymbol{x}(\cdot)$ satisfies the FDE in (29).

The inequality in (30) is the original definition of dissipativity established in [35]. In order to enforce the condition in (30) in conjunction with the Krasovskiĭ framework, however, we frequently utilize differential inequality

$$\widetilde{\forall} t \geq t_0, \quad \tfrac{\mathsf{d}}{\mathsf{d}t} \mathsf{v}(t, \mathbf{x}_t(\cdot)) - \mathsf{s}(\boldsymbol{z}(t), \boldsymbol{w}(t)) \leq 0 \qquad (31)$$

as a sufficient condition for (30) via an application of the Fundamental Theorem of Calculus for Lebesgue integrals, since $\dot{\mathsf{v}}(\mathbf{x}_t(\cdot))$ is integrable and exists for almost all $t \geq t_0$. Now the dissipativity inequality in (31) usually has to be addressed separately from the stability theorem if conventional derivatives are utilized for the differential equation. With Theorem 1, however, we can address (30) and the conditions in (5) together, as both of them are formulated using weak derivatives. Finally, it is worthy of mentioning that the procedures for utilizing Theorem 1 remain identical to those in the conventional Krasovskiĭ stability theorem, as evidenced in the derivations in [23].

## IV. Conclusion

A "measure" version of Krasovskiĭ stability theorem is established in this note for FDEs in the extended sense as defined in (1), satisfying the Carathéodory conditions. The theorem can be applied to FDEs for which the right-hand side mapping $\boldsymbol{f}(\cdot, \bullet)$ is not required to be continuous or piecewise continuous in $t$. Our theorem's proof features detailed explanations, highlighting key steps in (16)–(18), which eliminates the need to use the mean value theorem for differentiable functions. To show the advantage of the proposed theorem, we have presented examples of time-delay systems that can be addressed by Theorem 1 but were not possible using the conventional Krasovskiĭ stability theorem for the FDEs with ordinary derivatives.

## Acknowledgement

The authors sincerely thank Prof. Pierdomenico Pepe for introducing us to his excellent works in [28]–[30], [32].